\documentclass[reqno]{amsart}    
\usepackage{pstricks,pst-node,pst-coil,pst-plot}   

\theoremstyle{definition}      
\newtheorem{step}{Step} 
\newtheorem{thm}{Theorem}[section]     
\newtheorem{theorem}[thm]{Theorem}     
     
\newtheorem{corollary}[thm]{Corollary}     
     
\newtheorem{lemma}[thm]{Lemma}     
\newtheorem{prop}[thm]{Proposition}

\theoremstyle{remark}      
 
\newtheorem{remark}[thm]{Remark} 
     
\newtheorem{rem}[thm]{Remark}     
\theoremstyle{definition}      
     
\newtheorem{definition}[thm]{Definition}     
\def\al{{\alpha}}         
\def\be{{\beta}}         
\def\de{{\delta}}

\def\Om{{\Omega}}         
\def\la{{\lambda}}

\def\si{{\sigma}}         
         
\def\ga{{\gamma}}         
\def\ep{{\epsilon}}         
         
\def\th{{\theta}}         
\def\Th{{\Theta}}         
         
\def\phi{{\varphi}}

\let\pa\partial

\DeclareMathAlphabet{\doba}{U}{msb}{m}{n}

\gdef\mR{\doba{R}}         
\gdef\mS{\doba{S}}

\def\Spin{{\mathop{\rm Spin}}}

\def\Vol{{\mathop{\rm Vol}}}     
\let\vol\Vol
\def\Scal{{\mathop{\rm Scal}}}     
\let\scal\Scal

\def\Spec{{\mathop{\rm Spec}}}     

\def\End{{\mathop{\rm End}}}

\def\d2dt{\frac{d^2}{dt^2}}

\def\eref#1{{\rm (\ref{#1})}}   

\def\Rmax{R_{\textrm{max}}}

\let\<\langle 
\let\>\rangle

\long\def\ignorethis#1{}

\newdimen\templaenge

\def\Atbox#1#2{\setbox0\hbox{$\displaystyle #1$}\templaenge=\textwidth\advance\templaenge by -\wd0%
\setbox1\hbox{$#2$}\advance\templaenge by -\wd1%
$$#1\hbox{\kern\templaenge$#2$\hss}$$\par\bigbreak}

\newtheorem*{caseI}{Case I}
\newtheorem*{subcaseI.1}{Subcase I.1}
\newtheorem*{subcaseI.2}{Subcase I.2}
\newtheorem*{caseII}{Case II}
\newtheorem*{subsubcaseII.1.1}{Subsubcase II.1.1}
\newtheorem*{subsubcaseII.1.2}{Subsubcase II.1.2}
\newtheorem*{subcaseII.1}{Subcase II.1}
\newtheorem*{subcaseII.2}{Subcase II.2}
\newtheorem*{subsubcaseII.2.1}{Subsubcase II.2.1}
\newtheorem*{subsubcaseII.2.2}{Subsubcase II.2.2}
\newtheorem*{subsubcaseII.2.1.1}{Subsubcase II.2.1.1}
\newtheorem*{subsubcaseII.2.1.2}{Subsubcase II.2.1.2}
\def\an{a}
\def\detwo{\de_0}
\def\WS{$WS$}

\def\gWS{g_{\rm WS}}

\newcommand{\definedas}{\mathrel{\raise.095ex\hbox{\rm :}\mkern-5.2mu=}}
\begin{document}
\title{A surgery formula for the second Yamabe invariant}

\begin{abstract}
Let $(M,g)$ be a compact Riemannian manifold of dimension $n\geq 3$. For a metric $g$ on $M$, we let $\la_2(g)$ be the second eigenvalue of the Yamabe operator
$L_g:= \frac{4(n-1)}{n-2} \Delta_g +  \scal_g$. Then, the second Yamabe invariant is defined as 
$$
\si_2(M) \definedas \sup \inf_{h \in [g]}  \la_2(h) \Vol(M,h)^{2/n}. 
$$
where the supremum is taken over all metrics $g$ and the infimum is taken over the metrics in the conformal class $[g]$. Assume that  $\si_2(M)>0$. In the spirit of \cite{ammann.dahl.humbert:08}, we prove that if $N$ is obtained from $M$ by a $k$-dimensional surgery ($0 \leq k  \leq n-3$), there exists a positive constant 
$\Lambda_n$ depending only on $n$ such that 
$\si_2(N) \geq \min(\sigma_2(M), \Lambda_n)$. We then give some topological conclusions of this result.   
\end{abstract}
\maketitle     
\begin{center}
\sc S. El Sayed\footnote{elsayed@iecn.u-nancy.fr
}
\end{center}

\begin{center}
\today
\end{center}

\tableofcontents
\section{Introduction}
\noindent \textbf{Definition of the Yamabe operator $L_g$, eigenvalues of $L_g$, smooth Yamabe invariant $\sigma(M)$}\\
Let $(M,g) $ be a compact Riemannian manifold of dimension $n\geq 3$. We denote the scalar curvature by $\Scal_g$.   
Let us define 
$$\mu(M,g):= \inf_{\tilde{g} \in [g]} \int_M \Scal_{\tilde{g}} dv_{\tilde{g}} \left( \Vol_{\tilde{g}}(M) \right)^{-(n-2)/n}$$
and 
$$\sigma(M):= \sup_g \mu(M,g)$$
where, in the definition of $\mu(M,g)$,  the infimum runs over all the 
metrics $g'$ in the conformal class $[g]$ of $g$ and  
where, in the definition of $\sigma(M)$,  the supremum is taken over all the Riemannian metrics $g$ on $M$. The number $\mu(M,g)$,
also denoted by $\mu(g)$ if no ambiguity, is called the 
{\it Yamabe constant} while $\sigma(M)$ is called the {\it Yamabe invariant}. The Yamabe constant played a crucial role
 in the solution of the Yamabe problem 
solved between 1960 and 1984 by Yamabe, Tr\"udinger, Aubin and Schoen. This problem consists in finding a metric $\widetilde g$ conformal to 
$g$ such that the scalar curvature $\Scal_{\widetilde g}$ of $\widetilde g$ is constant. 
For more information, the reader  may refer to \cite{lee.parker:87, hebey:97, Aubin:98}. 
An important geometric meaning of $\mu(M,g)$ and $\sigma(M)$ is contained in the following well known result: 

\begin{prop}
Let $M$ be a compact differentiable manifold of dimension $n\geq 3$. Then, 
\begin{itemize} 
 \item if $g$ is a Riemannian metric on $M$, the conformal class $[g]$ of $g$ contains a metric of positive scalar curvature if and only if $\mu(M,g)>0$. 
\item $M$ carries a metric $g$ with positive scalar curvature if and only if $\sigma(M)>0$.
\end{itemize}
\end{prop} 
 
\noindent Classifying compact manifolds admitting a positive scalar curvature metric is a hard open problem which was studied by many mathematicians. 
Significant progresses were made thanks to surgery techniques.
We recall briefly that a surgery on $M$ is the procedure of constructing from $M$ a new manifold
$$N:= M\setminus {S^k\times B^{n-k}} \cup_{S^k\times S^{n-k-1}} \bar{B}^{k+1}\times S^{n-k-1},$$ 
by removing the interior of $S^k\times B^{n-k}$ and gluing it with  $\bar{B}^{k+1}\times S^{n-k-1}$ along the boundaries. 
Gromov-Lawson and Schoen-Yau proved in \cite{gromovlawson} and \cite{schoenyau} the following 
\begin{thm}
Let $M$ be a compact manifold of dimension $n \geq 3$ such that $\si(M) >0$.  
Assume that $N$ is obtained from $M$ by a surgery of dimension 
$k$ ($0 \leq k \leq n-3$). Then, $ \si(N) >0$. 
\end{thm}

\noindent Using cobordism techniques, one deduces: 
\begin{corollary}
Every manifold $M$ of dimension $n\geq 5$ simply connected and non-spin, carries a metric of positive scalar curvature.
\end{corollary}

\noindent Later, Kobayashi \cite{kobayashi} and Petean-Yun \cite{peteanyun} obtained new surgery formulas for $\sigma(M)$. These works were generalized by B. Ammann, M. Dahl and E. Humbert in \cite{ammann.dahl.humbert:08} where they proved in particular
\begin{thm} \label{adhl}
If $N$ is obtained from $M$ by a surgery of dimension $0\leq k \leq n-3$, then 
$$\sigma(N)\geq \min(\sigma (M), \Lambda_n),$$
where $\Lambda_n$ is a positive constant depending only on $n$.
\end{thm}
\noindent As a corollary, they obtained the following 
\begin{corollary} \label{coradh} 
Let $M$ be a simply connected compact manifold of dimension $n\geq 5$, then one of this assumptions is satisfied\\
\begin{enumerate}
\item $\sigma(M) = 0$ (which implies that $M$ is spin);
\item $\sigma (M)\geq \alpha_n$, where $\alpha_n$ is a positive constant depending only on $n$.
\end{enumerate} 
\end{corollary}

\noindent Now, let us define the \emph{Yamabe operator} or \emph{conformal Laplacian}  
$$L_g:= a\Delta_g + \Scal_g,$$
where $a = \frac{4(n-1)}{n-2}$ and where $\Delta_g$ is the Laplace-Beltrami operator. The operator $L_g$ is an elliptic differential operator of second order 
whose   spectrum is discrete:  
$$\Spec(L_g) = \{\la_1(g),\la_2(g),\cdots\},$$
where  $\la_1(g)< \la_2(g)\leq \cdots$ are the eigenvalues of $L_g$.  
The variational characterization of $\la_i(g)$ is given by 
$$\la_i(g) = \inf_{V\in Gr_i(H_1^2(M))}\sup_{v\in V\setminus\{0\}}\frac{\int_M vL_gv\, dv_g}{\int_M v^2\,dv_g},$$
where $Gr_i(H_1^2(M))$ stands for the $i$-th dimensional Grassmannian in $H_1^2(M).$ One important property of the eigenvalues of $L_g$ is that their sign is a conformal invariant equal to the sign of the Yamabe constant (see \cite{elsayed}). Consequently, 
a compact manifold $M$ possesses a metric with positive 
$\la_1$ if and only if it admits a positive scalar curvature metric. \\

\noindent Now, if $\mu(M,g)\geq 0,$ it is easy to check that 
\begin{eqnarray} \label{defmu} 
 \mu(M,g) = \inf_{\widetilde g\in \left[ g\right]} \la_1(\widetilde g) \Vol
(M,\widetilde g)^\frac{2}{n},
\end{eqnarray}
where $\left[ g\right]$ is the conformal class of $g$ and $\la_1$ is the first
eigenvalue of the Yamabe operator $L_g$. Inspired by these definitions,
one can define the \emph{second Yamabe constant} and the \emph{second Yamabe invariant} by 
$$\mu_2(M,g)=  \inf_{\widetilde g\in \left[ g\right]} \la_2(\widetilde g) \Vol
(M,\widetilde g)^\frac{2}{n},$$
and 
$$\si_2(M) = \sup_g \mu_2(M,g). $$

\noindent The second Yamabe constant $\mu_2(M,g)$ or $\mu_2(g)$ if no ambiguity was introduced  and studied in \cite{ammannhumbert} when $\mu(M,g) \geq 0$. This study was enlarged in \cite{elsayed} where we started to investigate the relationships between the sign of the second eigenvalue of the Yamabe operator $L_g$ and the existence of nodal
solutions of the equation $L_g u = \ep |u|^{N-2}u,$ where $\ep = -1, 0, +1$. The present paper establishes a surgery formula for $\si_2(M)$ in the spirit of Theorem \ref{adhl}. More precisely, our main result is the following
\begin{thm} \label{mainthm}
Let $M$ be a compact  manifold of dimension $n\geq 3$ such that  $\sigma_2(M)>0$. 
Assume that $N$ is obtained from $M$ by a surgery of dimension $0\leq k\leq n-3$, then we have
$$\sigma_2(N)\geq \min(\sigma_2(M), \Lambda_n),$$
where $\Lambda_n$ is a positive constant depending only on $n$.
\end{thm}

\noindent Note that B\"ar and Dahl in \cite{baerdahl} proved a surgery formula for the spectrum of the Yamabe operator with interesting topological consequences. \\

\noindent The proof of Theorem \ref{mainthm} is  inspired by the one of Theorem \ref{adhl} but some new difficulties arise here. Let us recall the strategy: first, we fix a metric $g$ on $M$ such that $\mu_2(M,g)$ is close to $\sigma_2(M)$. 
Then the goal is to construct on $N$ a sequence of metrics $g_\th$ such that 
$$ \liminf_{\th \to 0} \mu_2(N,g_\th) \geq \min(\mu_2(M,g),  \Lambda_n)$$ 
where $\Lambda_n >0$ depends only on $n$ (see Theorem \ref{theoremprincipal}). Surprisingly, if $\mu(M,g) =0$, we are not able to prove Theorem \ref{theoremprincipal} directly. So the first step is to show that one can assume that $\mu(M,g) \not= 0$ (see Paragraph \ref{munot0}).  
Here, we use exactly the same metrics than in \cite{ammann.dahl.humbert:08} and use many of their properties established in   \cite{ammann.dahl.humbert:08}. The proof consists in  studying the first and second eigenvalues $\la_1(u_\th^{N-2} g_\th)$ and $\la_2(u_\th^{N-2} g_\th)$ of $L_{u_ \th^{N-2} g_\th}$ where $u_\th$ is such that
$$\mu_2(g_\th) = \la_2(u_\th^{N-2} g_\th) \vol_{u_\th^{N-2} g_\th}(M)^{2/n},$$
or in other words, $u_\th$ is such that the metric $u_\th^{N-2} g_\th$ achieves the infimum in the definition of $\mu_2(N,g_\th)$. 
Two main difficulties arise in this situation: 
\begin{itemize} 
\item Contrary to what happened in \cite{ammann.dahl.humbert:08}, we could not show that  $\la_1(u_\th^{N-2} g_\th)$ and $\la_2(u_\th^{N-2} g_\th)$ are bounded.
\item The proof of Theorem \ref{adhl} was consisting in obtaining some good ``limit equations``. The difficulty here is to ensure that  
$$\lim_\th \la_1(u_\th^{N-2} g_\th) \not= \lim_\th \la_2(u_\th^{N-2} g_\th).$$ 
\end{itemize} 
The way to overcome these difficulties is to proceed in two steps: the first one is to show that  $ \la_2(u_\th^{N-2} g_\th)>0$. In a second step, we are able to get the desired inequality.\\

\noindent Let us now come back to Theorem \ref{mainthm}. Standard cobordism techniques  allow to deduce the following corollary
\begin{corollary}\label{cor}
Let $M$ be a compact, spin, connected and simply connected manifold of dimension $n\geq 5$ with  $n\equiv 0, 1, 2, 4$ mod $8$.
If $|\alpha(M)|\leq 1$, then 
$$\sigma_2(M)\geq \alpha_n,$$
where $\alpha_n$ is a positive constant depending only on $n$ and $\alpha(M)$ is the $\alpha$-genus of $M$ (see Section \ref{topologicalpart}).
\end{corollary}
\noindent  When $M$ is not spin, the conclusion of the corollary still holds but is a direct application of Corollary \ref{coradh} and the fact that $\sigma_2(M)\geq \sigma(M)$. Note that:\\

\noindent $\bullet$ In dimensions $1, 2$ mod $8$, $\alpha(M)\in \mathbb{Z}/{2\mathbb{Z}}$ and hence the condition on the $\alpha$-genus $|\alpha(M)|\leq 1$ is always satisfied. We then obtain that on any connected, simply connected manifold (not necessarily spin) of dimension $n\equiv 1, 2$ mod $8$
$$\sigma_2(M)\geq \alpha_n,$$
for some $\alpha_n>0$ depending only on $n$.\\
$\bullet$ In dimensions $0$ mod $8$, when $M$ is spin, $\alpha (M)= \hat{A}(M),$ where $\hat{A}$ is the $\hat{A}$-genus. Hence if $M$ is simply connected (not necessarily spin) connected of dimension $n\equiv 0$ mod $8$, $|\hat{A}|\leq 1$ then 
$$\sigma_2(M)\geq \alpha_n,$$
where $\alpha_n$ is a positive constant depending only on $n$.\\
$\bullet$ In dimensions $4$ mod $8$, when $M$ is spin, we have $\alpha(M) = \frac{1}{2}\hat{A}(M)$. When $M$ is spin and $\hat{A}(M)\leq 2$, we get that 
$|\alpha(M)|\leq 1$ and consequently,  for any   simply connected (not necessarily spin) connected $M$ of dimension $n\geq 5$, $n \equiv 4$ mod $8$
 with  $|\hat{A}|\leq 2$, we obtain that 
$$\sigma_2(M)\geq \alpha_n,$$ 
where $\alpha_n$ is a positive constant depending only on $n$.\\

\textbf{Acknowledgements:} I would like to thank Emmanuel Humbert for his encouragements, support and remarks along this work. I am also very grateful to Bernd Ammann, Mattias Dahl, Romain Gicquaud and Andreas Hermann for their remarks and their suggestions.


\section{Joining manifolds along a submanifold} \label{joining_man}
\subsection{Surgery on manifolds}
\begin{definition}
A surgery on a $n$-dimensional manifold $M$ is the procedure of constructing a new $n$-dimensional manifold 
$$N = (M\setminus f(S^k\times B^{n-k}))\cup (\overline{B}^{k+1}\times S^{n-k-1})/\sim,$$
by cutting out $f(S^k\times B^{n-k})\subset M$ and replacing it by $\overline{B}^{k+1}\times S^{n-k-1},$
where $f: S^k\times \overline {B^{n-k}}\rightarrow M$ is a smooth embedding which preserve the orientation and $\sim$ means that we paste along the boundary. Then, we construct on the topological space $N$ a differential structure and an orientation that makes a differentiable manifold such that the following inclusions 
$$M\setminus f(S^k\times B^{n-k})\subset N,$$
and 
$$\overline{B^{k+1}}\times S^{n-k-1}\subset N$$
preserve the orientation. We say that $N$ is 
obtained from $M$ by a surgery of dimension $k$ and we will 
denote $M \stackrel{k}{\rightarrow} N.$ 
\end{definition}
\noindent Surgery can be considered from another point of view. In 
fact, it is a special case of the connected sum: We paste $M$ and $S^n$ 
along a $k$-sphere. In this section we describe how two manifolds are joined
along a common submanifold with trivialized normal bundle. Strictly
speaking this is a differential topological construction, but since we
work with Riemannian manifolds we will make the construction adapted
to the Riemannian metrics and use distance neighborhoods defined by
the metrics etc. 
Let $(M_1,g_1)$ and $(M_2,g_2)$ be complete Riemannian manifolds of
dimension $n$. Let $W$ be a compact manifold of dimension $k$, where 
$0 \leq k \leq n$. Let $\bar{w}_i: W \times \mR^{n-k} \to TM_i$,
$i=1,2$, be smooth embeddings. We assume that $\bar{w}_i$ restricted
to $W \times \{ 0 \}$ maps to the zero section of $TM_i$ (which we
identify with $M_i$) and thus gives an embedding $W \to M_i$. The
image of this embedding is denoted by $W_i'$. Further we assume that
$\bar{w}_i$ restrict to linear isomorphisms $\{ p \} \times \mR^{n-k}
\to N_{\bar{w}_i(p,0)} W_i'$ for all $p \in W_i$, where $N W_i'$
denotes the normal bundle of $W_i'$ defined using $g_i$.
We set $w_i \definedas \exp^{g_i} \circ \bar{w}_i$. This gives
embeddings $w_i: W \times B^{n-k}(\Rmax) \to M_i$ for some 
$\Rmax > 0$ and $i=1,2$. We have $W_i' = w_i(W \times \{ 0 \})$ and we
define the disjoint union
$$ 
(M,g) \definedas (M_1 \amalg M_2, g_1 \amalg g_2),
$$
and 
$$
W' \definedas W_1' \amalg W_2'.
$$
Let $r_i$ be the function on $M_i$ giving the distance to $W_i'$.
Then $r_1 \circ w_1 (p,x) = r_2 \circ w_2(p,x) = |x|$ for $p \in W$,
$x \in B^{n-k}(\Rmax)$. Let $r$ be the function on $M$ defined by
$r(x) \definedas r_i(x)$ for $x \in M_i$, $i=1,2$. For $0 < \ep$ we
set $U_i(\ep) \definedas \{ x \in M_i \, : \, r_i(x) < \ep \}$ and
$U(\ep) \definedas U_1(\ep) \cup U_2(\ep)$. For $0 < \ep < \th$ we
define
$$
N_{\ep} 
\definedas
( M_1 \setminus U_1(\ep) ) \cup ( M_2 \setminus U_2(\ep) )/ {\sim},
$$
and 
$$
U^N_\ep (\th)
\definedas
(U(\th) \setminus U(\ep)) / {\sim}
$$
where ${\sim}$ indicates that we identify $x \in \partial U_1(\ep)$
with $w_2 \circ w_1^{-1} (x) \in \partial U_2(\ep)$. Hence
$$
N_{\ep} 
=
(M \setminus U(\th) ) \cup U^N_\ep (\th).
$$

\noindent We say that $N_\ep$ is obtained from $M_1$, $M_2$ (and $\bar{w}_1$, 
$\bar{w}_2$) by a connected sum along $W$ with parameter $\ep$. 

\noindent The diffeomorphism type of $N_\ep$ is independent of $\ep$, hence we
will usually write $N = N_\ep$. However, in situations when dropping
the index causes ambiguities, we will keep the notation $N_\ep$. For
example the function $r: M \to [0,\infty)$ gives a continuous function
$r_\ep: N_\ep \to [\ep, \infty)$ whose domain depends on $\ep$. It is
also going to be important to keep track of the subscript $\ep$ on 
$U^N_\ep (\th)$ since crucial estimates on solutions of the Yamabe
equation will be carried out on this set.  

\noindent The surgery operation on a manifold is a special case of taking
connected sum along a submanifold. Indeed, let $M$ be a compact
manifold of dimension $n$ and let $M_1 = M$, $M_2 = S^n$, 
$W = S^k$. Let $w_1 :  S^k \times B^{n-k} \to M$ be an embedding
defining a surgery and let $w_2 :  S^k \times B^{n-k} \to S^n$ be the
canonical embedding. Since $ S^n \setminus w_2 (S^k \times B^{n-k})$ is
diffeomorphic to $\overline{B^{k+1}} \times  S^{n-k-1}$ we have in this situation
that $N$ is obtained from $M$ using surgery on $w_1$, see 
\cite[Section VI, 9]{kosinski:93}.
\section{The constants $\Lambda_{n,k}$}\label{lambdank}
\subsection{Definition of $\Lambda_{n,k}$} In this paragraph, we define some constants $\Lambda_{n,k}$ in the same way than in \cite{ammann.dahl.humbert:08}. The only difference is that the  functions we considered are not necessarily positive. More precisely,
let $(M,h)$ be a Riemannian manifold of dimension $n\geq 3$. For $i = 1, 2$ we denote by $\Omega^{(i)}$ the set of $C^2$ functions $v$ (not necessarily positive) solution of the equation 
$$L_h v = \mu |v|^{N-2} v,$$ 
where  $\mu \in \mR$ . We assume that $v$ satisfies\\
\begin{eqnarray*}
&\bullet& v\not\equiv 0,\\
 &\bullet& \|v\|_{L^N(M)}\leq 1,\\
 &\bullet& v\in L^{\infty}(M),\\
\end{eqnarray*}
together with \\
\begin{eqnarray*}
&\bullet& v\in L^2(M), \text{ for } i = 1,\\
or\\
&\bullet& \mu\|v\|_{L^\infty(M)}^{N-2}\geq \frac{(n-k-2)^2(n-1)}{8(n-2)}, \text{ for } i = 2.
\end{eqnarray*}
For $i = 1, 2$, we set 
$$\mu^{(i)}(M,h) := \inf_{v\in \Omega^{(i)}(M,h)} \mu(v).$$\\
If $\Omega^{(i)}(M, h)$ is empty, we set $\mu^{(i)} = \infty.$
\begin{definition}
For $n\geq 3$ and $0\leq k \leq n-3$, we define 
$$\Lambda_{n,k}^{(i)} := \inf_{c\in [-1,1]}\mu^{(i)}({\mathbb H}_c^{k+1} \times {\mathbb S}^{n-k-1}),$$
and 
$$\Lambda_{n, k} := \min(\Lambda_{n,k}^{(1)}, \Lambda_{n, k}^{(2)}),$$
where $${\mathbb H}_c^{k+1} := ({\mathbb R}^k\times \mathbb R, \eta_c^{k+1} = e^{2ct}\xi^k + dt^2)$$
 \end{definition}
\noindent When considering only positive functions $v$, B. Ammann, M. Dahl and E. Humbert proved in  \cite{ammann.dahl.humbert:08} that these constants are positive. It is straightforward to see that the positivity of $v$ has no role in their proof and hence it remains true that $\Lambda_{n,k}>0$. 
 They gave also  explicit positive lower bounds of these constants and many of their techniques still hold in this context but we will not discuss this fact here.  For more informations, the reader may refer to \cite{ammanndahlhumbertlow},  \cite{ammanndahlhumbertsquare} and \cite{ammanndahlhumbetyamabeconstant} . 
\section{Limit spaces and limit solutions}
\begin{lemma}\label{vtheta}
Let $M$ be an $n$-dimensional manifold. let $(g_\th)$ be a sequence of metrics which converges toward a metric $g$ in $C^2$ on all compact $K\subset M$ when $\th\to 0$. Assume that $v_\th$ is a sequence of functions such that $\|v_\th\|_{L^\infty(M)}$ is bounded and 
$\|L_{g_\th} v_\th\|_{L^\infty(M)}$ tends to $0$. Then, there exists a smooth function $v$ solution of the equation 
$$L_g v = 0$$ 
such that $v_\th$ tends to $v$ in $C^1$ on each compact set $K \subset\subset V$. 
\end{lemma}

\noindent \textbf{Proof:} Let $K, K'$ be compact sets of $M$ such that $K'\subset K$, we have 
$$-g_\th^{ij}\left(\partial_i\partial_j v_\th-\Gamma_{ij}^k\partial_k v_\th\right)+\frac{n-2}{4(n-1)}\Scal_{g_\th} v_\th = f_\th\to 0.$$
Using Theorem 9.11 in \cite{gilbarg.trudinger:77}, one easily checks that
$$\|v_\th\|_{H^{2,p}(K',g)}\leq C(\|L_{g_\th} v_\th\|_{L^p(K,g_\th)} + \|v_\th\|_{L^p(K,g_\th)}).$$
It follows that $v_\th$ is bounded in $H^{2,p}(K',g)$ for all $p \geq 1$. Using  Kondrakov's theorem,  there exists $v_{K'}$ such that $v_\th$ tends to $v_{K'}$ in $C^1(K').$ Taking an increasing sequence of compact sets $K_m$ such that $\cup_m K_m = M$, $(v_\th)$ converges to $v_m$ on $C^1(K_m),$ we define $v:= v_m$ on $K_m$. Using the diagonal extraction process, we deduce that $v_\th$ tends to $v$ in $C^1$ on any compact set and that $v$ verifies the same Yamabe equation as $v_\th$. Since for each compactly supported smooth function $\phi$, we have 
$$\int_M L_{g_\th} \phi v_\th dv_{g_\th} \to \int_M L_g \phi v dv_g,$$ and $$\|L_{g_\th} v_\th \|_{L^\infty(M)} \to 0,$$  we obtain  that $L_g v = 0$ in the sense of distributions. Using standard regularity theorems, $v$ is smooth. 
\section{$L^2$-estimates on \WS-bundles} \label{wsbundles}
\noindent We suppose that the product $P \definedas I \times W \times S^{n-k-1}$ is equipped with a metric $\gWS$ of
the form 
$$
\gWS
=
dt^2 + e^{2\phi(t)}h_t + \sigma^{n-k-1}
$$
and we mean by \WS-bundle this product, where $h_t$ is a smooth family of metrics on $W$ and depending on $t$ 
and $\phi$ is a function on $I$. Let $\pi : P \to I$  be the projection onto the first factor and $F_t = \pi^{-1}(t) = \lbrace t\rbrace \times W\times S^{n-k-1}$, and the metric induced on $F_t$ is defined by 
$$g_t := dt^2 + e^{2\phi(t)} h_t + \sigma^{n-k-1}.$$
Let $H_t$ be the mean curvature of $F_t$ in $P$, it is given by the following 
$$H_t = -\frac{k}{n-1}\phi'(t)+ e(h_t),$$
with $e(h_t):= \frac{1}{2}tr_{h_t}(\pa_t h_t).$ The derivative of the element of volume of $F_t$ is 
$$\pa_t dv_{g_t} = -(n-1)H_tdv_{g_t}.$$
From the definition of $H_t$, when $t\to h_t$ is constant, we obtain that $$H_t = -\frac{k}{n-1}\phi'(t).$$
\begin{definition}
We say that the condition $(A_t)$ is verified if the following assumptions are satisfied:
\Atbox{\ \ 
\begin{matrix}
1.)& t \mapsto h_t\mbox{ is constant},\hfill\\ 
2.)& e^{-2 \phi(t)} \inf_{x \in W} \Scal^{h_t}(x) 
\geq -\frac{n-k-2}{32} \an ,\hfill\\
3.)& |\phi'(t)| \leq 1,\hfill\\
4.)& 0 \leq -2k \phi''(t) \leq \frac1{2} (n-1)(n-k-2)^2.
\end{matrix}}{(A_t)}
Similarly, for the condition $B_t$, we should have another assumptions to verify
\Atbox{\ \ 
\begin{matrix} 
1.) & t \mapsto \phi(t)\mbox{ is constant,}\hfill\\
2.) & \inf_{x\in F_t} \Scal^{\gWS}(x) 
\geq \frac{1}{2} \Scal^{\si^{n-k-1}}
= \frac{1}{2} (n-k-1)(n-k-2),\hfill\\
3.) & \frac{(n-1)^2}{2} e(h_t)^2 
+ \frac{n-1}{2} \partial_t e(h_t) 
\geq 
- \frac{3}{64} (n-k-2).\hfill
\end{matrix}}{(B_t)}
\end{definition}
\begin{theorem}\label{theo.fibest}
Let $\alpha,$ $\beta \in \mR$ such that $\left[ \alpha, \beta \right]\subset I.$ We suppose also that one of the conditions $(A_t)$ and $(B_t)$ is satisfied. We assume that we have a solution $v$ of the equation 
\begin{equation} \label{eqyamodif} 
L^{\gWS} v
= 
\an\Delta^{\gWS} v + \Scal^{\gWS} v
=
\mu u^{N-2}v + d^* A(dv) + Xv + \ep \pa_t v - sv
\end{equation}
where $s, \ep \in C^\infty(P)$, $A\in \End(T^*P)$, and $X\in \Gamma(TP)$ 
are perturbation terms coming from the difference between $G$ and
$\gWS$. We assume that the endomorphism $A$ is symmetric and that $X$
and $A$ are vertical, that is $dt(X) = 0$ and $A(dt) = 0$. Such that
\begin{equation} \label{u}
\mu \| u \|_{L^\infty(P)}^{N-2} 
\leq 
\frac{(n-k-2)^2(n-1)}{8(n-2)}.
\end{equation}
Then there exists $c_0>0$ independent of $\al$, $\beta$, and $\phi$,  
such that if 
$$
\| A \|_{L^\infty(P)}, 
\| X \|_{L^\infty(P)}, 
\| s \|_{L^\infty(P)},
\| \ep \|_{L^\infty(P)},
\| e(h_t) \|_{L^\infty(P)}
\leq 
c_0
$$
then 
$$
\int_{\pi^{-1} \left((\al + \ga ,\be - \ga)\right) }    
v^2 \, dv_{\gWS} 
\leq 
\frac{4 \| v \|_{L^\infty}^2}{n-k-2}  
\left(
\Vol^{g_\al} ( F_{\al }) + \Vol^{g_\be} ( F_{\be })
\right), 
$$
where $\ga \definedas \frac{\sqrt{32}}{n-k-2}$.
\end{theorem}
\noindent Remark that we should have $\beta - \alpha > 2\gamma$ to obtain our result and note that this theorem gives us an estimate of $\left\|v\right\|_{L^2}.$\\
For the proof of this Theorem, we mimic exactly the proof of Theorem 6.2 in \cite{ammann.dahl.humbert:08}. The only difference is that we consider here a nodal solution (and not a positive solution) of the equation $$L^{\gWS} v
=
\mu u^{N-2}v + d^* A(dv) + Xv + \ep \pa_t v - sv.$$
Other details are exactly the same.
\section{Main Theorem}
\noindent Theorem \ref{mainthm} is a direct corollary of 
\begin{theorem}\label{theoremprincipal}
Let $(M,g)$ be a compact Riemannian manifold  of dimension $n\geq 3$ such that $\mu_2(M,g)>0$ and let $N$ be obtained from $M$ by a surgery of dimension $0\leq k\leq n-3$. Then there exists a sequence of metrics $g_\th$ such that 
$$\liminf_{\th\to 0} \mu_2(N,g_\th)\geq \min (\mu_2(M,g), \Lambda_n),$$
where $\Lambda_n>0$ depends only on $n$.
\end{theorem}
\noindent Indeed, to get Theorem \ref{mainthm}, it suffices to apply Theorem \ref{theoremprincipal} with a metric $g$ such that 
$\mu_2(M,g)$ is arbitrary closed to $\sigma_2(M)$. The conclusion easily follows since $ \mu_2(N,g_\th) \leq \sigma_2(M)$. This section is devoted to the proof of Theorem \ref{theoremprincipal}.

\subsection{Construction of the metric $g_\th$}\label{constructionofthemetric}
\subsubsection{Modification of the metric $g$} \label{munot0} 
For a technical reason, we will need in the proof of Theorem \ref{theoremprincipal} that $\mu(g) \not=0$. To get rid of this difficulty, we need the following proposition:
\begin{prop} \label{munot0prop}
There exists on $M$ a metric $g'$ arbitrary close to $g$ in $C^2$ such that $\mu(g') \not= 0$. 
\end{prop}
Indeed, let us assume for a while that  Theorem \ref{theoremprincipal} is true if $\mu(g) \not=0$ and let us  
 see that the result remains true if $\mu(g) = 0$. A first observation is that if $g'$ is close enough to $g$ in $C^2$, then as one can check, $\mu_2(g')$ is close to $\mu_2(g)$.   Let us consider a metric $g'$ given by Proposition \ref{munot0prop} close enough to $g$ so that $\mu_2(g') > \mu_2(g) - \ep >0$ for an arbitrary small $\ep$.  From Theorem \ref{theoremprincipal} applied to  $g'$, we obtain a sequence of metrics $g_\th$  on $N$ such that 
$$\liminf_{\th\to 0} \mu_2(N,g_\th)\geq \min (\mu_2(M,g'), \Lambda_n) \geq \min (\mu_2(M,g) - \ep, \Lambda_n).$$
Letting $\ep$ tend to $0$, we obtain Theorem \ref{theoremprincipal}. 
It remains to prove Proposition \ref{munot0prop}. \\

\noindent {\bf Proof of Proposition \ref{munot0prop}: } At first, in order to simplify notations, we will consider $g$ as a metric on $M \amalg S^n$ and equal to the standard metric $g=\si^n$ on $S^n$. Since $\mu(g) = 0$, we can assume that $\scal_g = 0$, possibly making a conformal change of metrics. Let us consider a metric $h$ for which $\scal_h$ is negative and constant and whose existence is given in \cite{Aubin:98}.  Consider the analytic family of metrics $g_t:= t h + (1-t) g$. Since the first eigenvalue $\la_t$ of $L_{g_t}$ is simple, the function $t \to \la_t$ is analytic (see for instance Theorem VII.3.9 in \cite{kato:95}). Since $\la_0=0$ and $\la_1<0$, it follows that for $t$ arbitrary close to $0$, $\la_t \not= 0$.   Proposition \ref{munot0prop} follows since $\mu(g_t)$ has the same sign than $\la_t$.

\subsubsection{Definition of the metric $g_\th$}\label{definitionofthemetric}
\noindent As explained above, we will use the same construction as in \cite{ammann.dahl.humbert:08}. Consequently, we give the definition of $g_\th$ without additional explanations. The reader may refer to \cite{ammann.dahl.humbert:08} for more details. 
We keep the same notations than in Section \ref{joining_man}. Let $h_1$ be the restriction of $g$ to the surgery sphere $S_1'\subset M$ and $h_2$ be the restriction of the standard metric $\si^n=g$ on $S^n$ to $S_2' \subset S^n$. Define $S':= S_1' \amalg S_2'$ and 
$h \definedas h_1 \amalg h_2$ on $S'$. In the following, $r$ denotes the distance function to $S'$ in $(M\amalg S^n,g\amalg \si^n)$. In polar coordinates, the metric $g$ has the form 
\begin{equation} \label{metric=product} 
g =  h + \xi^{n-k} + T = h + dr^2 + r^2 \sigma^{n-k-1} + T
\end{equation}
on $U(\Rmax)\setminus S'\cong S'\times (0,\Rmax)\times
S^{n-k-1}$. Here $T$ is a symmetric $(2,0)$-tensor vanishing on $S'$ which is the error term measuring the fact that $g$ is not in general a product metric (at least near $S_1'$).
We also define the product metric
\begin{equation} \label{def.g'}
g' \definedas h + \xi^{n-k} = h + dr^2 + r^2 \si^{n-k-1},
\end{equation}
on $U(\Rmax) \setminus S'$ so that $g = g' + T$. As in \cite{ammann.dahl.humbert:08}, we have
\[ \left\{ \begin{array}{ccc} \label{normC0T}
| T(X,Y) | & \leq & Cr | X |_{g'} | Y |_{g'}, \\
|(\nabla_U T)(X,Y) | 
& \leq  & 
C | X |_{g'} | Y |_{g'} | U|_{g'},  \\ 
|(\nabla^2_{U,V}) T(X,Y) | 
& \leq& 
C | X |_{g'} | Y |_{g'} |U|_{g'}|V|_{g'}, 
\end{array} \right. \]
for $X,Y,U,V \in T_x M$ and  $x \in U(\Rmax)$. We define $T_1 \definedas T|_{M}$ and $T_2\definedas T|_{S^n}$. 
We fix $R_0\in (0,\Rmax)$, $R_0<1$ 
and  choose a smooth positive function $F: M \setminus S ' \to \mR$ such that 
$$
F(x) = 
\begin{cases}
1, 
&\text{if $x \in M \setminus U_1(\Rmax)\amalg S^n\setminus U_2(\Rmax)$;} \\ 
r(x)^{-1}, 
&\text{if $x \in U_i(R_0)\setminus S'$.}
\end{cases}
$$
Next we choose a sequence $\theta = \theta_j$ of positive numbers tending to $0$. For any $\theta$ we then choose a number  $\detwo = \detwo(\th) \in (0,\th)$ small enough to suit with the arguments below.
For any $\th>0$ and sufficiently small $\detwo$ there is 
$A_\th\in [\th^{-1}, (\detwo)^{-1})$ and a smooth function 
$f: U(\Rmax) \to \mR$ depending only on the coordinate $r$ such
that

$$
f(x)  =  
\begin{cases}
         -  \ln r(x),  &\text{if $x \in U(\Rmax) \setminus U(\th)$;} \\
\phantom{-} \ln A_\th, &\text{if $x \in U(\detwo)$,} 
\end{cases}
$$
and such that
\begin{equation} \label{asumpf}
\left| r\frac{df}{dr} \right|
= 
\left| \frac{df }{d(\ln r)} \right|
\leq 1, 
\quad 
\text{and} 
\quad
\left\|r\frac{d}{dr}\left(r\frac{df}{dr}\right)\right\|_{L^\infty}
= 
\left\|\frac{d^2f}{d^2(\ln r)}\right\|_{L^\infty}
\to 0
\end{equation}
as $\th\to 0$. 
Set $\ep = e^{-A_\th} \detwo$ that we assume smaller than $1$ and 
use this $\ep$ to construct $M$ as in Section \ref{joining_man}.
On $U^N_\ep(\Rmax) = 
\left( U(\Rmax) \setminus U(\ep) \right)/{\sim}$ we define $t$ by 
$$
t \definedas
\begin{cases}
         -  \ln r_1 + \ln \ep, & 
\text{on $U_1(\Rmax) \setminus U_1(\ep)$;} \\
\phantom{-} \ln r_2 - \ln \ep, & 
\text{on $U_2(\Rmax) \setminus U_2(\ep)$.}
\end{cases}
$$
One checks that  
\begin{itemize}
\item $
r_i = e^{|t|+ \ln \ep} = \ep e^{|t|};
$
\item 
$
F(x) = \ep^{-1}e^{-|t|}
$ 
for $x \in U(R_0) \setminus U^N(\th)$, or equivalently if 
$|t|+\ln \ep \leq \ln R_0$ and hence  
$$
F^2 g =\ep^{-2} e^{-2|t|}(h+T) + dt^2 + \sigma^{n-k-1}
$$
on $U(R_0)\setminus U^N(\th)$;
\item and 
$$
f(t) =  
\begin{cases}
-|t|-\ln\ep,  
&\text{if $\ln\th- \ln \ep \leq |t| \leq \ln\Rmax - \ln\ep$;} \\
\ln A_\th,  
&\text{if  $|t|\leq \ln\detwo-\ln\ep$.}
\end{cases}
$$
\end{itemize}
We have $|df/dt|\leq 1$, $\|d^2f/dt^2\|_{L^\infty}\to 0$. Now, we choose a
cut-off function $\chi:\mR\to [0,1]$ such that $\chi=0$ on
$(-\infty,-1]$, $|d\chi| \leq 1$, and $\chi=1$ on $[1,\infty)$. Finally, we define
$$
g_{\th}
\definedas  
\begin{cases}
F^2 g_i, 
&\text{on $M_i \setminus U_i(\th)$;} \\
e^{2f(t)}(h_i+T_i) + dt^2 + \sigma^{n-k-1},
&\text{on $U_i(\th)\setminus U_i(\detwo)$;} \\
\begin{aligned}
&A_{\th}^2 \chi( t / A_{\th} )(h_2+T_2)
+ A_{\th}^2 (1-\chi( t / A_{\th} ))(h_1+T_1)\\
&\quad + dt^2 + \sigma^{n-k-1}, 
\end{aligned}
&\text{on $U^N_\ep(\detwo)$.} 
\end{cases}
$$
Moreover, the metric $g_\th$ can be written as 
$$g_\th:= g'_\th + \widetilde {T_t} \text{ on } U^N(R_0),$$
where $g'_\th$ is the metric without error term and it is equal to 
$$g'_\th = e^{2f(t)}\widetilde{h_t} + dt^2 + \sigma^{n-k-1},$$
where the metric $\widetilde{h_t}$ is given by 
$$\widetilde{h_t}:= \chi(\frac{t}{A_\th})h_2 + (1-\chi(\frac{t}{A_\th}))h_1,$$
and $\widetilde{T_t}$ is the error term and his expression is given by the following 
$$\widetilde{T_t}:= e^{2f(t)}(\chi(\frac{t}{A_\th})T_2 + (1-\chi(\frac{t}{A_\th}))T_1).$$
We further have the following properties of the error term $\widetilde{T_t}$
\[ \left\{ \begin{array}{ccc} \label{norm}
|\widetilde T(X,Y) | & \leq & Cr | X |_{g'_\th} | Y |_{g'_\th}, \\
|\nabla {{\widetilde T}_t}|_{g'_\th} 
& \leq  & 
C e^{-f(t)},  \\ 
|\nabla^2{{\widetilde T}_t}|_{g'_\th} 
& \leq& 
C e^{-f(t)}, 
\end{array} \right. \]
where $\nabla$ is the Levi-Civita connection with respect to the metric $g'_\th$, for all $X$, $Y\in T_x N$ and $x\in U^N(R_0).$
\subsection{A preliminary result}
\noindent In order to prove Theorem \ref{theoremprincipal}, we will start by proving the following results. 
\begin{theorem}\label{theoremprincipal1}

\noindent {\bf Part 1:}  let $(u_\th)$ be a sequence of functions which satisfy
$$L_{g_\th}u_\th = \la_\th |u_\th|^{N-2} u_\th,$$
such that $\int_N |u_\th|^N dv_{g_\th}=1$ and $\la_\th\to_{\th\to 0} \la_\infty$, where $\la_\infty \in \mR$. Then, at least one of the two following assertions is true
\begin{enumerate}
 \item $\la_\infty\geq \Lambda_n$, where $\Lambda_n>0$ depends only on $n$;
\item there exists a function $u\in C^{\infty}(M \amalg S^n)$, $u \equiv 0$ on $S^n$, $u \not\equiv 0$ on $M$  solution of 
$$L_{g} u = \lambda_\infty |u|^{N-2} u,$$ 
with 
$$\int_M |u|^N dv_g=1$$
such that for all compact sets $K\subset M \amalg S^n \setminus S'$ (note that $K$ can also  be considered as a subset of $N$), $F^\frac{n-2}{2} u_\th$ tends to $u$ in $C^2(K)$, where $F$ is defined in Section \ref{constructionofthemetric}. Moreover, we have 
\begin{enumerate}
\item the norm $L^2$ of $u_\th$ is bounded uniformly in $\th$; \label{0}
\item $\lim_{b\to 0} \limsup_{\th\to 0} \sup_{U^N(b)} u_\th= 0$; \label{1}
\item $\lim_{b\to 0} \limsup_{\th\to 0}\int_{U^N(b)} u_\th^N\, dv_{g_\th} = 0.$ \label{2}
\end{enumerate}
\end{enumerate}
\noindent {\bf Part 2:} let $u_\th$ be as in Part 1 above and assume that Assertion 2) is true. Let $v_\th$ be a sequence of functions which satisfy 
$$L_{g_\th} v_\th = \mu_\th |u_\th|^{N-2} v_\th,$$
such that $\int_N v_\th^N dv_{g_\th} = 1,$ $\mu_\th\to \mu_\infty$ where $\mu_\infty < \mu(\mS^n)$. Then, there exists a function $v \in C^{\infty}(M \amalg S^n)$, $v \equiv 0$ on $S^n$, $v \not\equiv 0$ on $M$ solution of 
$$L_g v = \mu_\infty |u|^{N-2} v$$
with 
$$\int_M |v|^N dv_g = 1$$
and such that for all compact sets $K\subset M \amalg S^n \setminus S' $, $F^{\frac{n-2}{2}} v_\th$ tends to $v$ in $C^2(K).$ Moreover,
\begin{enumerate}
\item  the norm $L^2$ of $v_\th$ is bounded uniformly in $\th$; \label{2'}
\item $\lim_{b\to 0} \limsup_{\th\to 0} \sup_{U^N(b)} v_\th= 0;$ \label{3}
\item $\lim_{b\to 0} \limsup_{\th\to 0} \int_{U^N(b)} v_\th^N\, dv_{g_\th} = 0.$ \label{4}
\end{enumerate}
\end{theorem}
\subsubsection{Proof of Theorem \ref{theoremprincipal1} Part 1} \label{part1}
Let $(u_\th)$ be a sequence of functions which satisfy
$$L_{g_\th}u_\th = \la_\th |u_\th|^{N-2} u_\th,$$
such that $\int_N |u_\th|^N \,dv_{g_\th}=1$ and $\la_\th\to_{\th\to 0} \la_\infty$, where $\la_\infty \in \mR$.
We proceed exactly as in \cite{ammann.dahl.humbert:08} where here, the manifold $M_2$
is $S^n$ equiped with the standard metric $\sigma^n$, and where $W$ is the sphere $S^k$. The only difference will be that $u_\th$ may now have a changing sign.

\begin{remark}
 In the proof of the main theorem in \cite{ammann.dahl.humbert:08},  it was proven that  
$$\la_\infty > -\infty.$$ 
 Here, we made the assumption that $\la_\infty$ has a limit. Without this assumption, one could again prove that 
$\la_\infty > -\infty$ but the point here is that there is no reason why $\la_\infty$ should be bounded from above contrary to what happened in   \cite{ammann.dahl.humbert:08}. 
\end{remark}

\noindent The argument of Corollary 7.7 in \cite{ammann.dahl.humbert:08} still holds here and shows that 
\begin{eqnarray} 
 \liminf_\th \|u_\th\|_{L^\infty(N)} > 0. 
\end{eqnarray}
Several cases are studied: 
\begin{caseI}
$\limsup_{\th\to 0}\|u_\th\|_{L^\infty(N)} = \infty$.
\end{caseI}
\noindent Set $m_\th:= \|u_\th\|_{L^\infty(N)}$ and choose $x_\th\in N$ such that $u_\th(x_\th) = m_\th$. After taking a subsequence, we can assume that $\lim_{\th\to 0} m_\th = \infty$. We have to study the following two subcases.
\begin{subcaseI.1}
There exists $b>0$ such that $x_\th\in N \setminus U^N(b)$ 
for an infinite number of $\th$.
\end{subcaseI.1}
\begin{subcaseI.2}
For all $b>0$ it holds that $x_\th\in U^N(b)$ for $\th$ sufficiently
small.
\end{subcaseI.2}
\begin{caseII}
There exists a constant $C_0$ such that 
$\| u_\th\|_{L^{\infty}(N)}\leq C_0$ for all $\th$. 
\end{caseII}
\begin{subcaseII.1}
 There exists $b>0$ such that
$$
\liminf_{\th\to 0} 
\left( \la_\th\sup_{ U^N(b) } {u_\th}^{N-2} \right)
< 
\frac{(n-k-2)^2(n-1)}{8(n-2)}.
$$
\end{subcaseII.1}

\begin{subsubcaseII.1.1}
$\limsup_{b\to 0} \limsup_{\th\to 0} \sup_{U^N(b)} u_\th> 0$.  
\end{subsubcaseII.1.1}

\begin{subsubcaseII.1.2}
$\lim_{b\to 0} \limsup_{\th\to 0} \sup_{U^N(b)} u_\th= 0$.  
\end{subsubcaseII.1.2}

\begin{subcaseII.2}
$$
\la_\th  \sup_{ U^N(b) } {u_\th}^{N-2}  
\geq
\frac{(n-k-2)^2(n-1)}{8(n-2)}
$$
\end{subcaseII.2}
\noindent In Subcases I.1 and I.2, it is shown in \cite{ammann.dahl.humbert:08} that $\lambda_\infty \geq \mu(\mS^n)$. The proof still holds when $u_\th$ has a changing sign. In Subsubcase II.1.1 and Subcase 
II.2, we obtain that  $\lambda_\infty \geq \Lambda_{n,k}$ where $\Lambda_{n,k}$ is a positive number depending only on $n$ and $k$. The definition of $\Lambda_{n,k}$ in \cite{ammann.dahl.humbert:08} is the infimum of energies of positive solutions of the Yamabe equation on model spaces (see Section \ref{lambdank}). This definition has to be slightly modified to allow nodal solutions. As explained in Section \ref{lambdank} the proof that $\Lambda_{n,k}> 0$ remains the same.   \\


\noindent In Subcases I.1, I.2, II.1.1 and II.2, we then get that 
$\lambda_\infty\geq \Lambda_n$, where $$\Lambda_n:= \min_{k\in\{0, \cdots, n-3\}} \{\Lambda_{n,k}, \mu\}.$$ In particular, Assertion 1) of part 1 in Theorem \ref{theoremprincipal1} is true. So let us examine Subsubcase  II.1.2. The assumption of Subcase II.1 allows to obtain as in \cite{ammann.dahl.humbert:08} that 
\begin{eqnarray} \label{uboundedl2}
\int_N u_\th^2 dv_{g_\th} \leq C. 
\end{eqnarray} 
for some $C>0$. The  assumptions of Subcase II.1.2 are that 
\begin{eqnarray} \label{ubounded}
 \sup_N (u_\th) \leq C
\end{eqnarray}
and that 
\begin{eqnarray} \label{uto0}
\limsup_{b\to 0} \limsup_{\th\to 0} \sup_{U^N(b)} u_\th= 0. 
\end{eqnarray}
\begin{step} \label{step1} 
We prove that
$\lim_{b \to 0} \limsup_{\th\to 0} \int_{U^N(b)} {|u_\th|}^N \, dv_{g_\th}
= 0$. 
\end{step}
\noindent Let $b >0$. We have, by Relation \eref{uboundedl2} 
$$
\int_{U^N(b) } |u_\th|^N \, dv_{g_\th} 
\leq 
A_0  \sup_{U^N(b)} |u_\th|^{N-2}, 
$$  
where $A_0$ is a positive  number which does not depend on $b$ and $\th.$ The claim then follows from \eref{uto0}.
\begin{step}\label{step2}
$C^2$ convergence on all compact sets of $M\amalg S^n\setminus S'.$
\end{step}
\noindent Let $(\Omega_j)_j$ be an
increasing sequence of subdomains of $(M \amalg S^n\setminus S')$ with smooth boundary such that $\bigcup_{j}\Omega_j = M\amalg S^n\setminus S',$ $\Omega_j\subset \Omega_{j+1}$. The norm $\left\|u_\th\right\|_{L^\infty(N)}$ is bounded, 
then  so is $\left\|u_\th\right\|_{L^\infty(\Omega_{j+1})}$. Using  standard results on elliptic regularity (for more details, see for example \cite{gilbarg.trudinger:77}), we see that the sequence $(u_\th)$ is bounded in the Sobolev space 
$H^{2,p}(\Omega'_{j})$ $\forall p \in (1, \infty)$ where $\Omega'_j$ is any domain such that $\overline{\Om}_j \subset \Om'_j \subset \overline{\Om'_j} \subset \Om_{j+1}$.  The Sobolev embedding Theorem implies that $(u_\th)$ is bounded in $C^{1,\alpha}(\overline{\Omega_j})$ for any $\alpha \in (0, 1)$. (See Theorem 4.12 in \cite{adams.fournier:03} for more informations on Sobolev embedding Theorems).\\
Now we use a diagonal extraction process, by taking successive subsequences, it follows that $(u_\th)$ converges to functions $\widetilde{u_j}\in C^1(\overline{\Omega_j})$ and such that $\widetilde{u_j}\vert_{\overline{\Omega}_{j-1}} = \widetilde {u}_{j-1}.$\\
We define 
$$\widetilde u =\widetilde {u_j} \text{ on } \overline{\Omega_j}.$$ 
By taking a diagonal subsequence of $u_\th$, we get that $u_\th$ tends to $\widetilde u$ in $C^1$ on any compact subset of $M\amalg S^n\setminus S'$ and by $C^1$-convergence of the functions $u_\th$, the function $\widetilde u$ satisfies the equation 
\begin{equation}\label{eqV1}
L_{g_\th} \widetilde u =\la_\infty |\widetilde u|^{N-2}{\widetilde {u}} \text{ on } M\amalg S^n\setminus S'.
\end{equation}
We recall that  $g_\th= F^2 g = (F^\frac{n-2}{2})^\frac{4}{n-2}g$ on $U^N(b)$. 
By conformal invariance of the Yamabe operator we obtain for all $v$
$$L_{F^2 g}v = F^{-\frac{n+2}{2}}L_g(F^\frac{n-2}{2}v).$$
Now we set 
$$ u = F^\frac{n-2}{2} \widetilde u.$$
We obtain 
\begin{eqnarray*}
L_g u&=&F^\frac{n+2}{2}L_{F^2 g}\widetilde u\\
&=& F^\frac{n+2}{2} \la_{\infty} |\widetilde u|^{N-2}\widetilde u\\
&=& \la_{\infty} |u|^{N-2} u.
\end{eqnarray*}
This shows that $u$ is a solution on $(M\amalg S^n\setminus S', g)$ of the following equation
$$L_g u = \la_{\infty} |u|^{N-2} u.$$
Moreover, using Step \ref{step1} and the fact that $\int_N  u_\th^N dv_{g_\th}=1$, the function $u$ satisfies 
\begin{eqnarray*}
\int_{M\amalg S^n} u^N\, dv_g &=&\int_{M\amalg S^n \setminus S'} {\widetilde u}^N\, dv_g\\
&=& \lim_{b\to 0}\lim_{\th \to 0}\int_{U^N(b)} u_{\th}^N\, dv_{g_{\th}}\\
&=& 1.
\end{eqnarray*}
\begin{step}\label{step3}
Removal of the singularity
\end{step}
\noindent The next step is to show that $u$ is a solution on all $M\amalg S^n$ of 
\begin{eqnarray}\label{eqnwidetilde}
L_g {u} = \la_\infty |u|^{N-2}u.
\end{eqnarray}
To prove this fact, we will show that for all $\phi \in C^\infty(M\amalg S^n)$, we have
$$\int_{M\amalg S^n} L_g u \phi\,dv_g =\int_{M\amalg S^n}\la_\infty |u|^{N-2} u \phi\,dv_g.$$
First, we have 
\begin{eqnarray*}
\int_{M\amalg S^n} u L_g\phi\,dv_g &=&\int_{M\amalg S^n} u L_g(\phi-\chi_\ep \phi+\chi_\ep \phi)\,dv_g\\
&=&\int_{M\amalg S^n} u L_g(\chi_\ep \phi)\,dv_g+\int_{M\amalg S^n} u L_g((1-\chi_\ep)\phi)\,dv_g,
\end{eqnarray*}
where 
\begin{eqnarray*}
   \left|\;  
   \begin{matrix}
     \chi_\ep=1\hfill & \hbox{if } d_g(x,S')<\ep,\\\\
     \chi_\ep =0\hfill & \hbox{if } d_g(x,S')\geq 2\ep,\\\\
     \left| d\chi_\ep\right|<\frac{2}{\ep}.
   \end{matrix}
   \right.
\end{eqnarray*}
Since $(1-\chi_\ep)$ is compactly supported in $M\amalg S^n\setminus S'$, we have \begin{eqnarray*}
\int_{M\amalg S^n} u L_g((1-\chi_\ep)\phi)\,dv_g &=& \int_{M\amalg S^n} (L_g u)(1-\chi_\ep)\phi\,dv_g\\
&\to&
 \int_{M\amalg S^n} L_g u\phi\, dv_g = \int_{M\amalg S^n} \la_\infty |u|^{N-2} u\phi\,dv_g.
 \end{eqnarray*}
Then, it remains to prove that 
$$\int_{M\amalg S^n} u L_g(\chi_\ep \phi)\, dv_g\rightarrow 0.$$
We have 
\begin{eqnarray*}
L_g(\chi_\ep\phi)&=&C_n \Delta (\chi_\ep\phi)+\Scal_g(\chi_\ep\phi)\\
&=&C_n \Delta \chi_\ep \phi+ C_n\Delta \phi \chi_\ep+\Scal_g(\chi_\ep \phi)-2\left\langle \nabla \chi_\ep,\nabla \phi\right\rangle\\
&=& \chi_\ep L_g\phi+C_n(\Delta \chi_\ep)\phi-2\left\langle \nabla \chi_\ep,\nabla \phi\right\rangle.
\end{eqnarray*}
According to Lebesgue Theorem, it holds  that
$$\int_{M\amalg S^n} u \chi_\ep L_g\phi \, dv_g\rightarrow 0 \text{ a.e.}$$
Further, we have
\begin{eqnarray} \label{etoile} 
\left| \int_{M\amalg S^n} u L_g(\chi_\ep \phi)\,dv_g\right|&\leq& \frac{C}{\ep^2}\int_{C_\ep} u\,dv_g\\
&\leq& \frac{C}{\ep^2}\left(\int_{C_\ep} u^2\,dv_g\right)^\frac{1}{2}\left(\vol(Supp(C_\ep))\right)^\frac{1}{2},
\end{eqnarray}
where $C_\ep = \left\lbrace x\in M\amalg S^n; \ep < d(x,S')<2\ep \right\rbrace = U^N(2\ep)\setminus U^N(\ep)$.\\
In addition, we get from \eref{uboundedl2} that 

$$\int_N {\widetilde u}^2\,dv_{F^2 g} < +\infty,$$
which implies that
$$\int_{C_\ep}{\widetilde u}^2\,dv_{F^2 g} < +\infty.$$
Let us compute
\begin{eqnarray*}
\int_{C_\ep}{\widetilde u}^2\,dv_{g_{\th}}&=& \int_{C_\ep} \left(F^{\frac{n-2}{2}}\right)^{\frac{2n}{n-2}}F^{-(n-2)}u^2\,dv_g \\
&=& \int_{C_\ep} F^2 u^2\,dv_g < +\infty.
\end{eqnarray*}
We recall that $F = \frac{1}{r}$ on $C_\ep$. Coming back to \eref{etoile}, we deduce 
\begin{eqnarray*}
\left| \int_M u L_g(\chi_\ep \phi)\,dv_g\right|&\leq&\frac{C}{\ep^2} \left(\int_{C_\ep}\frac{u^2 F^2}{F^2} \,dv_g\right)^\frac{1}{2} \left(\vol(C_\ep)\right)^\frac{1}{2}\\
&\leq& \frac{C}{\ep^2}\times \ep \times \ep^\frac{n-k}{2} = C\ep^{\frac{n-k}{2}-1}.
\end{eqnarray*}
Since $k\leq n-3,$ we have 
$$\frac{n-k}{2}-1 > 0,$$
which implies that 
$$\int_{M\amalg S^n} u L_g(\chi_\ep \phi)\,dv_g \rightarrow 0.$$
Finally, we get that $u$ is a solution on $M\amalg S^n$ of the equation 
$$L_g u= \la_\infty |u|^{N-2} u.$$
\begin{step}
We have either $u \equiv 0$ on ${\mS}^n$ either $\la_\infty\geq \mu({\mS}^n).$
\end{step}
\noindent Note that the function $u$ verifies 
\begin{eqnarray}\label{mamalgs}
\int_{M\amalg S^n} |u|^N\, dv_g \leq 1.
\end{eqnarray}
Since
\begin{eqnarray*}
\int_{M\amalg S^n} |u|^N\, dv_g &=& \int_{M\amalg S^n} |\widetilde u|^N\, dv_{g_{\th}}\\
&\leq& \int_N |\widetilde u|^N \, dv_{g_{\th}}\\
&\leq& \lim_{\th\to 0} \int_N |{u_\th}|^N \,dv_{g_\th} = 1. 
\end{eqnarray*}
Assume that  $u \not\equiv 0$ on ${\mS}^n$.\\
Setting $w = u_{|{\mS^n}}$ and using equations (\ref{eqnwidetilde}) and (\ref{mamalgs}), we have 
\begin{eqnarray*}
\mu({\mS}^n) \leq Y(w)&=& \frac{\la_\infty \int_{{\mS}^n} w^N\, dv_g}{\left(\int_{{\mS}^n} w^N\, dv_g\right)^\frac{n-2}{n}}\\
&=& \la_\infty \left(\int_{{\mS}^n} w^N \,dv_g\right)^\frac{2}{n}\leq \la_\infty.
\end{eqnarray*}
Then we obtain that $\la_\infty \geq \mu({\mS}^n)$ and hence, the conclusion $1)$ of Theorem \ref{theoremprincipal1} Part 1 is true.

\subsubsection{Proof of Theorem \ref{theoremprincipal1} Part 2}
 We consider a function $v_\th$ satisfying 
\begin{eqnarray}\label{vth}
L_{g_\th} v_\th = \mu_\th |u_\th|^{N-2} v_\th,
\end{eqnarray}
with 
$$\int_N |v_\th|^N\, dv_{g_\th} = 1.$$
\noindent A first remark is the following: as in Lemma 7.6 of \cite{ammann.dahl.humbert:08}, we observe that $U^N(b)$ is a $WS$-bundle  for any $b > 0$. Since $u_\th$ satisfies 
$$\lim_{b\to 0}\limsup_{\th\to 0} \sup_{U^N(b)} u_\th = 0.$$
Then, for $b$ small enough, we have 
$$\mu_\th \|u_\th\|_{U^N(b)}^{N-2} \leq \frac{(n-k-2)^2(n-1)}{8(n-2)}.$$
We then can apply  Theorem \ref{theo.fibest} on $U^N(b)$ and the proof of  Lemma 7.6 of \cite{ammann.dahl.humbert:08} shows that there exists numbers $c_1, c_2>0$ independent of $\th$ such that 
 \begin{eqnarray}\label{c1c2}
\int_N |v_\th|^2 \, dv_{g_\th} \leq c_1{\Arrowvert v_\th\Arrowvert}^2_{L^\infty(N)} + c_2.
\end{eqnarray}
As a consequence, we get that 
$$
\liminf_{\th\to 0} \| v_\th\|_{L^{\infty}(N)} 
>0.
$$
Indeed, assume that 
$$\lim_{\th \to 0}\left\|v_\th\right\|_{L^\infty(N)} = 0.$$
By Equation \eref{c1c2}, we have
\begin{eqnarray*}
1 = \int_N |v_\th|^N \, dv_{g_\th}&\leq& {\Arrowvert v_\th\Arrowvert}_{L^\infty(N)}^{N-2} \int_N |v_\th|^2\, dv_{g_\th}\\
&\leq& {\Arrowvert v_\th\Arrowvert}_{L^\infty(N)}^{N-2} (c_1 {\Arrowvert v_\th\Arrowvert}_{L^\infty(N)}^2+ c_2)\to 0,
\end{eqnarray*}
as $\th\to 0$. This gives the desired contradiction.
 In the rest of the proof, we will study several cases. In what follows, only Subcase II.1.2 will be a big deal: Subcases I.1, I.2 and II.1 will be excluded by arguments mostly contained in \cite{ammann.dahl.humbert:08}. So we will just give few explanations for these cases.
 \begin{caseI}
$\limsup_{\th\to 0}\|v_\th\|_{L^\infty(N)} = \infty$.
\end{caseI}

\noindent Set $m_\th\definedas \|v_\th\|_{L^\infty(N)}$ and choose  $x_\th\in N$ with $v_\th(x_\th) = m_\th$. After taking a
subsequence we can assume that $\lim_{\th\to 0} m_\th= \infty$. 
\begin{subcaseI.1}
There exists $b>0$ such that $x_\th\in N \setminus U^N(b)$ 
for an infinite number of $\th$.
\end{subcaseI.1}
 \noindent By taking a subsequence we can assume
that there exists $\bar{x} \in M \amalg S^n \setminus U(b)$ such
that $\lim_{\th\to 0} x_\th= \bar{x}$. We define 
$\tilde{g}_\th\definedas m_\th^{\frac{4}{n-2}} g_\th$. For $r>0$, \cite{ammann.dahl.humbert:08} tells that for $\th$ small enough, there exists a diffeomorphism 
$$
\Th_\th: 
B^n(0,r)
\to
B^{g_\th} ( x_\th, m_\th^{-\frac{2}{n-2}} r) 
$$ 
such that the sequence of metrics 
$(\Th_\th^* (\tilde{g}_\th))$ tends to the flat 
metric $\xi^n$ in $C^2(B^n(0,r))$, where $B^n(0,r)$ is the standard ball in $\mR^n$ centered in $0$ with radius $r$. We let 
$\tilde{u}_\th\definedas m_\th^{-1} u_\th,$ ${\tilde{v}}_\th\definedas m_\th^{-1} v_\th$ and we have
\begin{eqnarray*}
L_{\tilde{g}_\th} {\tilde{v}}_\th &=& \la_\th{\tilde{u}}_{\th}^{N-2} {\tilde{v}}_\th\\
&=& \frac{\la_\th}{m_\th^{N-2}}{u_\th}^{N-2}\tilde{v}_\th.
\end{eqnarray*} 

\noindent Since $\left\|u_\th\right\|_{L^\infty(N)} \leq C,$ it follows that
$\left\|L_{\tilde{g}_\th} \tilde{v}_\th \right\|_{L^\infty(N)}$ tends to $0$.
Applying Lemma \ref{vtheta}, we obtain a solution $v\not\equiv 0$ of the following equation on $\mR^n$
$L_{\xi^n} v = 0$. Since $\Scal_{\xi^n} = 0$, $v$ is harmonic and admits a maximum at $x = 0$. As a consequence, $v$ is constant equal to $v(0) = 1$. This is a contradiction, since $\|v\|_{L^N}\leq 1.$
\begin{subcaseI.2}
For all $b>0$ it holds that $x_\th\in U^N(b)$ for $\th$ sufficiently
small.
\end{subcaseI.2}
\noindent We proceed as in Subcase I.2 in \cite{ammann.dahl.humbert:08}. As in Subcase I.1 above, we get from Lemma \ref{vtheta} a function $v$ which is harmonic on $\mR^n$ and admits a maximum at $x = 0$. This is again a contradiction.
\begin{caseII}
There exists a constant $C_0$ such that 
$\| v_\th\|_{L^{\infty}(N)}\leq C_0$ for all $\th$. 
\end{caseII}
\noindent By (\ref{c1c2}), there exists a constant $A_0$ independent of $\th$ such that  
\begin{equation} \label{vboundedl2}
\| v_\th\|_{L^2(N,g_\th)} \leq A_0.
\end{equation}
We split the treatment of  Case II into two subcases.
\begin{subcaseII.1}
$\limsup_{b\to 0} \limsup_{\th\to 0} \sup_{U^N(b)} v_\th> 0$.  
\end{subcaseII.1}
\noindent Again mimicking what is done in \cite{ammann.dahl.humbert:08}, we obtain from Lemma \ref{vtheta} a function $v$ which is a solution of 
$L_{G_c} v= 0$ on $\mR^{k+1} \times S^{n-k-1}, G_c$ for some $c \in [-1,1]$ where 
$G_c= e^{2cs} \xi^k + ds^2 + \sigma^{n-k-1}$. In Subcases I.1 and I.2, we used the fact that 
$\frac{\la_\th}{m_\th^{N-2}}$ tends to $0$ to show that at the limit $L_{G_c} v= 0$. Here, the argument is different: first we set $\alpha_0:= \frac{1}{2}\limsup_{b\to 0}\limsup_{\th\to 0} \sup_{U^N(b)} v_\th >0$. Then, we can suppose that there exists a sequence of positive numbers $(b_i)$ and $(\th_i)$ such that 
$$\sup_{U^N(b_i)}v_{\th_i}\geq \alpha_0,$$
for all $i$. To simplify, we write $\th$ for $\th_i$ and $b$ for $b_i$. Take $x'_\th\in \overline{U^N(b_\th)}$ such that 
$$v_\th(x'_\th)\geq \alpha_0.$$
For $r,r'>0$, we define 
$$U_\th(r,r'):= B^{{\widetilde h}_{t_\th}}(y_\th, e^{-f(t_\th)}r)\times [t_\th-r',t_\th+r']\times S^{n-k-1}.$$
As in \cite{ammann.dahl.humbert:08}, the function $v$ is obtained as the limit of $v_\th$ on each $U_\th(r,r')$ (with $r, r'>0$).
The fact that $L_{G_c} v = 0$ follows from the observation that 
$$\sup_{U_\th(r,r')} |u_\th| = 0,$$
hence 
$$|u_\th|^{N-2}v_\th \to 0 \text{ uniformly on } U_\th(r,r').$$

\begin{subcaseII.2}
$\lim_{b\to 0} \limsup_{\th\to 0} \sup_{U^N(b)} v_\th = 0$.  
\end{subcaseII.2}
\noindent By the same method than in Subsection {\ref{part1}}, we obtain that there is a function $v$ solution of the following equation
$$L_g v = \mu_\infty |u|^{N-2} v,$$
such that 
$$\int_N v^N\,dv_g\leq 1.$$
Suppose that $v\not\equiv 0$ on $\mS^n$, then we have
\begin{eqnarray*}
\mu(\mS^n) \leq Y(v) = \mu_\infty\frac{\int_{\mS^n} u^{N-2} v^2\, dv_g}{(\int_{\mS^n} v^N\, dv_g)^\frac{2}{N}} = 0
\end{eqnarray*}
since $u \equiv 0$ on $S^n$. This is a contradiction.
This proves that $v \not\equiv 0$ on $\mS^n$. By the same argument than in Part 1, we have 
$\int_M |v|^N dv_g = 1$.  We finally obtain that the function $v$ satisfies all the desired conclusions of Theorem \ref{theoremprincipal1} Part 2.

\subsection{Proof of Theorem \ref{theoremprincipal}}

\noindent Let $(g_\th)$ the sequence of metrics defined on $N$ as in Section {\ref{constructionofthemetric}}.\\
{\bf{Step 1:}} For $\th$ small enough, we show that if $$\la_k(M,g)>0 \Rightarrow \la_k(N,g_\th)>0,$$ where $\la_k$ is the $k^{th}$ eigenvalue associated to the Yamabe equation.\\

\begin{remark}
Note that this step implies that the existence of a metric with positive $\la_k$ is preserved by surgery of dimension $k \in \{0,\cdots,n-3\}$. This is an alternative proof of a result already contained in 
\cite{baerdahl}. 
\end{remark}

\noindent We proceed by contradiction and we suppose that $\la_k(N, g_\th)\leq 0.$ Let $u_\th$ be a minimizing solution of the Yamabe problem. By referring to \cite{elsayed}, there exists functions $v_{\th, 1}= u_\th, v_{\th,2}, \cdots, v_{\th, k}$ solution of the following equation on $N$
$$L_{g_{\th}} v_{\th, i} =  \la_{\th, i} u_\th^{N-2} v_{\th,i},$$
where $$\la_{\th, i} = \la_i(N, u_\th^{N-2} g_\th),$$
such that 
$$\int_N {v_{\th, i}}^N\, dv_{g_\th} = 1 \text{ and }\int_N {u_\th}^{N-2} v_{\th, i} v_{\th, j} \,dv_{g_\th}= 0 \text{ for all }i \neq j.$$
By conformal invariance of the sign of the eigenvalues of the Yamabe operator (see \cite{elsayed}), we have $$\la_{\th,i}=\la_i(N, u_\th^{N-2} g_\th)\leq 0.$$
Moreover, by construction, it is easy to check that $ \la_{\th,1}  = \mu_\th$ where $\mu_\th=\mu(N,g_\th)$ is the Yamabe constant of the metric $g_\th$. The main theorem in \cite{ammann.dahl.humbert:08} implies that  $\lim_{\th\to 0} \la_{\th, 1} = \lim_{\th\to 0} \mu_\th > -\infty$.  It follows that there exists a constant $C >0$ such that   $-C \leq \la_{\th, 1} \leq \cdots\leq\la_{\th, k} \leq 0$. Then, for all $i$, $\la_{\th, i}$ is bounded and by restricting to a subsequence we can assume that $\la_{\infty, i} := \lim_{\th\to 0} \la_{\th, i}$ exists. Parts 1) and 2) of Theorem \ref{theoremprincipal} give the existence of functions $u=v_1, \cdots, v_k$ defined on $M,$ with $v_i\neq 0$ for all $i$ such that $F^{\frac{n-2}{2}}v_{\th,i}$ tends to $v_i$ in $C^1$  on each compact set $K \subset M \amalg S^n \setminus S'$.  The functions $v_i$ are solutions of the following equation
$$L_g v_i = \la_{\infty, i} u^{N-2} v_i.$$
Moreover, we have
$$\int_M |v_i|^N\,dv_g\leq 1 \text{ and } \lim_{b\to 0}\limsup_{\th\to 0}\int_{U^N_\ep(b)}|v_{\th, i}|^N\, dv_g = 0.$$ 
Let us show that for all $i\neq j$, we get that  
$$\int_M u^{N-2} v_i v_j \,dv_g = 0.$$
Set 
$$\widetilde {u}_\th = F^\frac{n-2}{2} u_\th,$$
and 
$$\widetilde {v}_{\th, i} = F^\frac{n-2}{2} v_{\th, i}.$$
For $b>0$ small, we have for $i \not= j$
\begin{eqnarray*} 
\int_{M\setminus {U(b)}}u^{N-2} v_i v_j\, dv_g &  = \lim_{\th\to 0}\int_{M\setminus U(b) = N\setminus U^N_\ep(b)} \widetilde u_\th^{N-2} \widetilde v_{\th, i} \widetilde v_{\th, j} dv_{g} \\
&  =  \lim_{\th\to 0}\int_{M\setminus U(b) = N\setminus U^N_\ep(b)} u_\th^{N-2} v_{\th, i}  v_{\th, j} dv_{g_\th}
\end{eqnarray*}
where we used $dv_{g_\th} = F^n dv_g$. Using now the fact that $\int_N u_\th^{N-2} v_{\th, i} v_{\th, j}\, dv_{g_\th} = 0$,   we get
\begin{eqnarray*}
\left|\int_{M\setminus U(b)} u^{N-2} v_i v_j\, dv_g \right|&=& \left|\lim_{\th\to 0} \int_{N\setminus {U^N_\ep(b)}} u_{\th}^{N-2} v_{\th,i} v_{\th,j} \, dv_{g_\th}\right|\\
&=&\lim_{\th\to 0} \left|\int_{U^N_\ep(b)}u_\th^{N-2} v_{\th,i} v_{\th,j} \, dv_{g_\th}\right|.
\end{eqnarray*}
We write 
\begin{eqnarray*}
\left|\int_{U^N_\ep(b)} u_\th^{N-2} v_{\th,i} v_{\th,j} \, dv_{g_\th}\right|&\leq& \left(\int_{U_\ep^N(b)} u_{\th}^N\, dv_{g_\th}\right)^{\frac{N-2}{N}} \left(\int_{U_\ep^N(b)} |v_{\th,i}|^N\, dv_{g_\th}\right)^{\frac{1}{N}}\\
&&\left(\int_{U_\ep^N(b)} |v_{\th,j}|^N\, dv_{g_\th}\right)^{\frac{1}{N}}.
\end{eqnarray*}
Using the assertion 
$$\lim_{b\to 0}\limsup_{\th\to 0}\int_{U^N_\th(b)} v_{\th, i}^N\, dv_{g_\th} = 0.$$
we obtain that 
$$\lim_{b\to 0} \limsup_{\th\to 0} \left|\int_{U^N_\ep(b)} u_\th^{N-2} v_{\th, i} v_{\th, j}\, dv_{g_\th}\right| = 0.$$
We get finally that
$$\left|\int_M u^{N-2} v_i v_j\, dv_g\right| = \lim_{b\to 0}\left|\int_{M\setminus U(b)}u^{N-2} v_i v_j\, dv_g \right| = 0 \text{ for all }i\neq j.$$
We now write
\begin{eqnarray*}
0<\la_k(M,g)&\leq& \sup_{(\alpha_1,\cdots ,\alpha_k)\neq (0,\cdots ,0)}F(u, \alpha_1 v_1+\cdots +\alpha_k v_k)\\
&=& \sup_{(\alpha_1,\cdots ,\alpha_k)\neq (0, \cdots, 0)} \frac{\int_M (\alpha_1 v_1+\cdots +\alpha_k v_k)L_g(\alpha_1 v_1+\cdots +\alpha_k v_k)\,dv_g}{\int_M {u}^{N-2}(\alpha_1 v_1+\cdots +\alpha_k v_k)^2\,dv_g}\\
&=& \sup_{(\alpha_1,\cdots ,\alpha_k)\neq (0, \cdots, 0)} \frac{\alpha_1^2 \int_M v_1 L_g v_1\, dv_g + \cdots + \alpha_k^2 \int_M v_k L_g v_k\, dv_g}{\alpha_1^2\int_M u^{N-2} v_1^2\, dv_g+\cdots+\alpha_k^2\int_M u^{N-2}v_k^2\, dv_g}\\
&=&\sup_{(\alpha_1,\cdots ,\alpha_k)\neq (0, \cdots, 0)} \frac{\alpha_1^2\la_{\infty,1}\int_M u^{N-2} v_1^2\,dv_g +\cdots+\alpha_k^2\la_{\infty,k}\int_M u^{N-2}v_k^2\,dv_g}{\alpha_1^2\int_M u^{N-2} v_1^2\, dv_g +\cdots+\alpha_k^2\int_M u^{N-2} v_k^2\, dv_g}\\
&\leq& 0,
\end{eqnarray*}
since each $\la_{\infty,i}\leq 0$. This gives the desired contradiction.\\

\begin{remark} 
Note that, for $i \geq 2$ it could happen that $\int_M u^{N-2} v_i^2 dv_g = 0$ if $M$ is not connected.
\end{remark}

\noindent {\bf{Step 2:} Conclusion} \\
Since $\mu_2(M,g)>0,$ from Step 1, we get that $\mu_2(N, g_\th) > 0.$ Assume $\mu_2(N, g_\th) < \mu(\mS^n)$ (otherwise, we  are  done). Using \cite{elsayed} we construct a sequence $(v_\th)$ solution of 
$$L_{g_\th}v_\th = \mu_2(N,g_\th) |v_\th|^{N-2}v_\th,$$
such that 
$$\int_N v_\th^N\, dv_{g_\th} = 1.$$
 By Theorem \ref{theoremprincipal1} Part 1), this holds that $\lim_{\th \to 0}  \mu_2(N,g_\th) \geq \Lambda_n$ (and the conclusion of Theorem \ref{theoremprincipal} is true) or there exists a function $v$ solution on $M$ of the equation:
$$L_g v = \mu_\infty |v|^{N-2}v,$$
with $\mu_\infty = \lim_\th \mu_2(N,g_\th) \geq 0$ and 
$$\int_M |v|^N\, dv_g =  1.$$ 
This is what we assume until now. 

\noindent As explained in Paragraph \ref{munot0}, we can assume that $\mu(g) \not= 0$.\\ 

\noindent {\bf Case 1:} $\mu(g)<0$.
 
\noindent Assume that $M$ is connected (so is $N$)  and let us prove that $v$ has a changing sign. 
We suppose by contradiction that $v\geq 0.$ The maximum principle gives that $v>0$. Let $u$ be a positive solution of the Yamabe equation on $M,$ i.e.
$$L_g u = \mu(g) u^{N-1}.$$
Since $v>0$, we can write: 
$$L_g v = \underbrace{\mu_\infty}_{\geq  0} |v|^{N-2} v = \mu_\infty v^{N-1}.$$
Multiplying the second equation by $u$ and integrating, we get 
$$\underbrace{\mu(g)}_{<0}\int_M u^{N-1} v \,dv_g= \int_M L_g u v \, dv_g = \int_M uL_gv\,dv_g = \underbrace{\mu_\infty}_{\geq 0} \int_M v^{N-1}u\,dv_g.$$
This gives a contradiction.
Then $v$ have a changing sign and this implies that 
$$\mu_2(M,g)\leq \sup_{\alpha, \beta}F(v, \alpha v^+ + \beta v^-) = \mu_\infty.$$

 \noindent If $M$ is now disconnected, then the Yamabe minimizer $u$ is positive on a connected component of $M$. If $uv\not\equiv 0$, the same proof holds. If $uv \equiv 0$ then 
$$\mu_2(M,g)\leq \sup_{\alpha, \beta}F(v, \alpha u + \beta v) = \mu_\infty$$
 In any case, the conclusion of  Theorem \ref{theoremprincipal} is true.  \\

\noindent {\bf{Case 2:}} $\mu(M,g) >0$.

\noindent Then,  $\la_1(N,g_\th)>0$. In \cite{elsayed}, it is established that the sign of the eigenvalues of the Yamabe operator is conformally invariant. Consequently, 
$\la_1(N, v_\th^{N-2} g_\th)>0$. Set $\mu_1 = \la_1(N, v_\th^{N-2} g_\th)$ and let $u_\th$ be associated to $\mu_1$. Since associated to the first eigenvalue of the Yamabe operator, $u_\th$ is positive on at least one connected component of $N$ (and $0$ on the other). 	In addition, $u_\th$ is a solution  of the equation 
$$L_{g_\th}u_\th = \mu_1 |v_\th|^{N-2}u_\th,$$
such that 
$$\int_N u_\th^N\, dv_{g_\th} = 1 \text{ and } \int_N |v_\th|^{N-2} u_\th v_\th\,dv_{g_\th} = 0.$$
Using Theorem \ref{theoremprincipal1} Step 2), there exists a function $u$ solution on $M$ of the following equation
$$L_g u = \mu_{\infty, 1} |v|^{N-2} u,$$
where $\mu_{\infty,1} := \lim_\th \mu_1$. Note that this limit exists after a possible extraction of a subsequence since 
$0 \leq \mu_1 \leq \mu_2(N,g_\th)$. 
Proceeding as in Step 1, we show that 
\begin{eqnarray} \label{uvorth}
 \int_M |v|^{N-2}u v \, dv_g = 0.
\end{eqnarray}

\noindent By maximum principle and since $u_\th>0$,   $u>0$ on at least one connected component of $M$. 
Then, $u$ and $v$ satisfy the equations
$$L_g u = \mu_{\infty, 1} |v|^{N-2} u,$$
and 
$$L_g v = \mu_\infty |v|^{N-2}v.$$
These equations implies that $\mu_{\infty, 1}$ and $\mu_\infty$ are some eigenvalues of the generalized metric $|v|^{N-2} g$ (see \cite{elsayed}). Since positive, $u$ is associated to the first eigenvalue of $L_{|v|^{N-2} g}$ i.e. $\mu_{\infty, 1} = \la_1(M, |v|^{N-2} g)$. Hence, $\mu_{\infty, 1} \leq \mu_\infty$.\\
Finally, we obtain that
\begin{eqnarray*}
\mu_2(M,g ) &\leq&   \la_2(|v|^{N-2}g) \Vol_{|v|^{N-2}g}(M)^{\frac{2}{n}}= \mu_\infty
\end{eqnarray*}
since 
$$\Vol_{|v|^{N-2}g}(M) = \int_M |v|^N dv_g = 1$$
and since $\mu_{\infty, 1} \leq \mu_\infty$  are associated to two non proportional eigenfunctions in the metric $|v|^{N-2} g$ (thanks to Relation \eref{uvorth})
where we recall that $\mu_\infty = \lim_{\th\to 0} \mu_2(N, g_\th).$ This proves Theorem \ref{theoremprincipal}.
\begin{rem}
The reason why we need $\mu(g) \not= 0$ is the following. If $\mu(g) = 0$, the proof of Case 1 clearly does not lead to a contradiction. So, we would like to apply the method used in Case 2 above. For this,  we need that $\la_1(v_\th^{N-2} g_\th)$ is bounded. When $\mu(g) >0$, this holds true since  
$$0 \leq \la_1(v_\th^{N-2} g_\th) \leq \la_2(v_\th^{N-2} g_\th) = \mu_2(N,g_\th) \to \mu_\infty.$$ 
If $\mu(g) = 0$, one cannot say nothing about the sign of $\la_1(v_\th^{N-2} g_\th)$. In particular, if it is negative, we were not able to prove that $\la_1(v_\th^{N-2} g_\th)$ is bounded from above and the proof breaks down.   
\end{rem}

\section{Some applications}\label{topologicalpart}
\noindent In this section, we establish some topological applications of Theorem \ref{mainthm}.
\subsection{A preliminary result} 
\noindent We have 
\begin{prop}
Let $V$, $M$ be two compact manifolds such that $V$ carries a metric $g$ with $\Scal_g = 0$ and $\sigma (M)>0,$ then 
$$\sigma_2(V\amalg M)  \geq \min(\mu_2(g),\sigma(M))>0.$$
\end{prop}
\noindent \textbf{Proof:} On $V\amalg M$, let $G = \la g + \mu h$, where $\la$ and $\mu$ are two positive constants and for a small $\ep$, $h$ is a metric such that $\sigma(M) \leq\mu(M,h) + \ep$. We have
\begin{eqnarray*}
\Spec(L_G) &=& \Spec(L_{\la g}) \cup \Spec(L_{\mu h})\\
&=& \la^{-1}\Spec(L_g) \cup \mu^{-1}\Spec(L_h)\\
&=& \{\la^{-1} \la_1, \la^{-1}\la_2, \cdots \}\cup \{\mu^{-1} \la_1',\mu^{-1} \la_2', \cdots \}
\end{eqnarray*}
where $\la_i$ (resp. $\la_i'$) denotes the $i$-th eigenvalue of $L_g$ (resp. $L_h$). 
The assumption we made allows to claim that $\la_1=0$, $\la_2>0$ and $\la_1' >0$. Hence, we deduce that $\la_2(L_G) = \min \{\la^{-1}\la_2, \mu^{-1} \la_1'\}$.\\
We know that 
$$\vol_G(V\amalg M) = \la^\frac{n}{2} \vol_g(V) + \mu^\frac{n}{2} \vol_h(M).$$
$\bullet$ For $\mu = 1$ and $\la\to +\infty$, we have 
$$\la_2(L_G) = \la^{-1}\la_2.$$
\begin{eqnarray*}
\la_2(L_G) {\vol_G}^\frac{2}{n}(V\amalg M) &=& \la^{-1} \la_2 \left(C + \la {\vol_g}^\frac{2}{n}(V)\right)\\
&\to_{\la \to +\infty} & \la_2 {\vol_g}^\frac{2}{n}(V) = \mu_2(g).
\end{eqnarray*}
$\bullet$ For $\la = 1$ and $\mu\to +\infty$, in this case $$\la_2(L_G)) = \mu^{-1} \la_1'.$$
Hence \begin{eqnarray*}
\la_2(L_G) {\vol_G}^\frac{2}{n}_g(V\amalg M) &=& \mu^{-1} \la_1'\left(C + \mu {\vol_h}^\frac{2}{n}\right)\\
&\to_{\mu \to +\infty} & \la_1' {\vol_h}^\frac{2}{n} = \mu(M,h) \geq \sigma(M) - \ep.
\end{eqnarray*}
Finally we get that
$$\sigma_2(V\amalg M)\geq \min(\mu_2(g), \sigma(M)).$$
\begin{rem}
\begin{enumerate}
\item It is known that if $\sigma(M)>0$ and $\sigma(N)>0$, then 
$$\sigma(M\amalg N) = \min(\sigma(M), \sigma(N)),$$
where $M\amalg N$ is the disjoint union of $M$ and $N$. (see \cite{ammann.dahl.humbert:08}).\\ 
\item Let $V$ with $\sigma(V) \leq  0$, then for $k\geq 2$
$$\sigma_2(\underbrace{V\amalg\cdots\amalg V}_{k \hbox{ times }} \amalg M) \leq 0.$$
\end{enumerate}
\noindent Indeed, let any metric $g = g_1\amalg g_2\amalg \cdots \amalg g_k\amalg g_n$ on $V\amalg\cdots\amalg V \amalg M$. Let $v_i$ be functions associated to $\la_1(g_i)$ which is non-negative by assumption. The functions $\tilde {v_i} = 0 \amalg \cdots 0 \amalg \underbrace{v_i}_{i^{th} \hbox{factor} } 
\amalg\hskip0.1cm  0 \cdots \amalg 0$ are linearly independent and satisfy $L_g(\tilde{v_i}) =  \la_1(g_i) v_i$ and thus are eigenfunctions of $L_g$.  This implies that $\la_k(g) \leq 0$ and since $k \geq 2$, $\la_2(g) \leq 0$. 
\end{rem}
\noindent This remark explains the condition $|\al(M)| \leq 1$ in Corollary \ref{cor}: it is used to ensure that $M$ is obtained from a model manifold $V\amalg N$ with a number of factors $V$ (where $V$ carries a scalar flat metric and $\sigma(N)>0$) not larger than 1. We recall that  the $\alpha$-genus is an homomorphism from the spin cobordism ring $\Omega_*^{\Spin}$to the real $K$-theory ring $KO_*(pt),$
$$\alpha: \Omega_*^{\Spin} \rightarrow KO_*(pt).$$
It is important that $\alpha$ is a ring homomorphism, i.e. for any connected closed spin manifolds $M$ and $N$, $\alpha(M\amalg N) = \alpha(M) + \alpha(N)$ and $\alpha(M\times N) = \alpha(M).\alpha(N).$\\
Noting that $KO_n(pt)$ vanishes if $n = 3, 5, 6, 7$ mod $8$, is isomorphic to $\mathbb{Z}$ if $n = 0,4$ mod $8$ and is isomorphic to $\mathbb{Z}/{2\mathbb{Z}}$ if $n = 1, 2$ mod $8$. Recall also that $\alpha$ is exactly the $\hat{A}$-genus in dimensions $0$-mod $8$ and equal to $\frac{1}{2}\hat{A}$-genus in dimensions $4$ mod $8$.\\
In \cite{baerdahl1}, Proposition 3.5 says that in dimensions $n = 0, 1, 2, 4$ mod $8$, there exists a manifold $V$ such that $\alpha(V) = 1$ and $V$ carries a metric $g$ such that $\Scal_g = 0$.\\
$\bullet$ When $\alpha(M) = 0$ then Thm A in \cite{stolz} applies and $\sigma(M)\geq \alpha_n$ where $\alpha_n$ depending only on $n$.
\begin{thm}
Let $M$ be a spin manifold, if $\alpha(M) = 0$, this is equivalent to the existence of a manifold $N$ cobordant to $M$ such that the scalar curvature of $N$, $\Scal_g$ is positive. 
\end{thm}
\noindent Remember that a cobordism is a manifold $W$ with boundary whose boundary is partitioned in two, $W = M\amalg (-N)$.
\begin{thm}
If $M$ is cobordant to $N$ and if $M$ is connected then $M$ is obtained from $N$ by a finite number of surgeries of dimension $0\leq k\leq n-3.$
\end{thm}
\begin{prop}
Let $M$ be a spin, simply connected, connected manifold of dimension $n\geq 5$, if $n = 0, 1 ,2 , 4$ mod $8$ and $|\alpha(M)|\leq 1$, then
$$\sigma_2(M)\geq \alpha_n,$$
where $\alpha_n$ is a positive constant depending only on $n$. 
\end{prop}
\noindent \textbf{Proof:} Proposition 3.5 in \cite{baerdahl1} gives us that for each $n = 0, 1, 2, 4$ mod $8$, $n\geq 1$, there is a manifold $V$ of dimension $n$ such that $V$ carries a metric $g$ such that $\Scal_g = 0$ and $\alpha(V) = 1$.\\
$\bullet$ \textbf{First case:} If $\alpha(M) = 0$, then $M$ is cobordant to a manifold $N$ such that $\Scal_g$ on $N$ is positive. In  this case we can obtained $M$ from $N$ by a finite number of surgeries of dimension $k\leq n-3$. Hence, by Corollary  $\sigma(M)\geq c_n$ with $c_n$ is a positive constant depending only on $n$.\\
$\bullet$ \textbf{Second case:} If $\alpha(M) = 1$, then $\alpha(M\amalg (-V)) = 0$, so there exists a manifold $N$ with $\Scal_g > 0$ such that $M\amalg (-V)$ is cobordant to $N$ which is equivalent to say that $M$ is cobordant to $V\amalg N$. Consequently $M$ can be obtained from $V\amalg N$ by a finite number of surgeries of dimension $k\leq n-3$. Applying the main theorem of this paper, we get the desired result.


\begin{thebibliography}{10}
\bibitem{adams.fournier:03}
R.A.~Adams and J.J.F.~Fournier, \emph{Sobolev spaces}, 2nd edition, Pure and Applied Mathematics, 140, Academic Press, Amsterdam, 2003.

\bibitem{ammanndahlhumbertlow}
 B.~Ammann, M.~Dahl, and E.~Humbert, \emph{Low-Dimensional surgery and the Yamabe invariant}, arXiv:1204.1197.
 
\bibitem{ammanndahlhumbertsquare}
\bysame, \emph{Square-Integrability of solutions of the Yamabe equation}, arXiv:1111.2780.
  
\bibitem{ammann.dahl.humbert:08}
\bysame, \emph{Smooth Yamabe invariant and surgery}, Pr{\'e}publications de l'institut Elie Cartan, 58 pages, 2008.


\bibitem{ammanndahlhumbetyamabeconstant}
\bysame, \emph{The conformal Yamabe constant of product manifolds}, to appear in Proc. AMS.

\bibitem{ammannhumbert}
B.~Ammann and E.~Humbert, \emph{The second Yamabe invariant},  J. Funct. Anal., 235 (2006), No 2, p. 377-412.

  
  \bibitem{Aubin:98}
  T.~Aubin, \emph{Some nonlinear problems in {R}iemannian geometry}, Springer Monographs in Mathematics, Springer-Verlag, Berlin, 1998.
  
  \bibitem{baerdahl}
  C.~B{\"a}r and M.~Dahl, \emph{Small eigenvalues of the conformal Laplacian}, Journal: Geometric and Functional Analysis GAFA - Volume 13, Issue 3, p. 483-508.
  
 \bibitem{baerdahl1}
  \bysame, \emph{Surgery and the spectrum of the Dirac operator}, J. Reine Angew. Math., \textbf{552} (2002), 53--76.
  
  
  \bibitem{elsayed}
  S.~El Sayed, \emph{Second eigenvalue of the Yamabe operator and applications}, arXiv:1204.1268.
  
\bibitem{gilbarg.trudinger:77}
D.~Gilbarg and N.~Trudinger, \emph{Elliptic partial differential equations of
  second order}, Grundlehren der mathematischen Wissenschaften, no. 224,
  Springer-Verlag, 1977.
  
\bibitem{gromovlawson}
M.~Gromov and H. B.~Lawson, \emph{The classification of simply connected manifolds of positive scalar curvature}, Ann. of Math. (2) 111 (1980), 423--434, MR0577131, Zbl 0463.53025.
  

\bibitem{hebey:97}
E.~Hebey, \emph{{Introduction \`a l'analyse non-lin\'eaire sur les
  vari\'et\'es}}, Diderot Editeur, Arts et Sciences, out of print, no longer
  sold, 1997.
  
\bibitem{kato:95}
T.~Kato, \emph{{Perturbation theory of linear operators}}, Reprint of the 1980
edition, Springer Verlag, Berlin, Heidelberg (1995).

\bibitem{kobayashi}
O.~Kobayashi, \emph{Scalar curvature of a metric with unit volume}, Math. Ann. \textbf{279} (1987), 253--265, MR0919505, Zbl 0611.53037.
  
 \bibitem{kosinski:93}
A.~A. Kosinski, \emph{Differential manifolds}, Pure and Applied Mathematics,
  vol. 138, Academic Press Inc., Boston, MA, 1993.
 
 \bibitem{lee.parker:87}
J.~M. Lee and T.~H. Parker, \emph{{The Yamabe problem}}, Bull. Am. Math. Soc.,
  New Ser. \textbf{17} (1987), 37--91.

 \bibitem{peteanyun}
  J.~Petean and G.~Yun, \emph{Surgery and the Yamabe invariant}, Geom. Funct. Anal. \textbf{9} (1999), 1189--1199, MR1736933, Zbl 0976.53045.
  
  \bibitem{schoenyau}
  R.~Schoen and S.-T.~Yau, \emph{On the structure of manifolds with positive scalar curvature}, Manuscripta Math. \textbf{28} (1979), 159--183, MR0535700. Zbl 0423.53032.
  
  \bibitem{stolz}
S.~Stolz, \emph{Simply connected manifolds of positive scalar curvature}, Ann. Math., II. Ser. \textbf{36} (1992), 511--540.
  
\end{thebibliography}
\end{document}